\documentclass[11pt]{article}
\usepackage{amsmath, amssymb, amscd, amsthm, amsfonts}
\usepackage{graphicx}
\usepackage{hyperref}
\usepackage{enumitem}
\usepackage{tikz}

\title{Counterexamples to a conjecture of Gao and Gu on weighted Catalan numbers}
\author{Ryota Inagaki \and Dimana Pramatarova}
\date{}

\newtheorem{theorem}{Theorem}
\newtheorem{lemma}[theorem]{Lemma}
\newtheorem{definition}[theorem]{Definition}
\newtheorem{conjecture}[theorem]{Conjecture}

\begin{document}

\maketitle
\renewcommand{\thefootnote}{\fnsymbol{footnote}} 



\begin{abstract}
We construct infinitely many counterexamples to a conjecture of Gao and Gu on when the 2-adic valuations of weighted Catalan numbers equal those of Catalan numbers.
\end{abstract}

\section{Introduction}
\subsection{Background}
The sequence of Catalan numbers $\{C_n\}_{n\geq 0}$ is one of the most famous integer sequences in combinatorics. The $n$th Catalan number is given by the formula $C_n = \frac{1}{n+1}\binom{2n}{n}$. Catalan numbers count Dyck paths, triangulations of convex $(n+2)$-gons, and many other objects \cite{MR3467982}. 

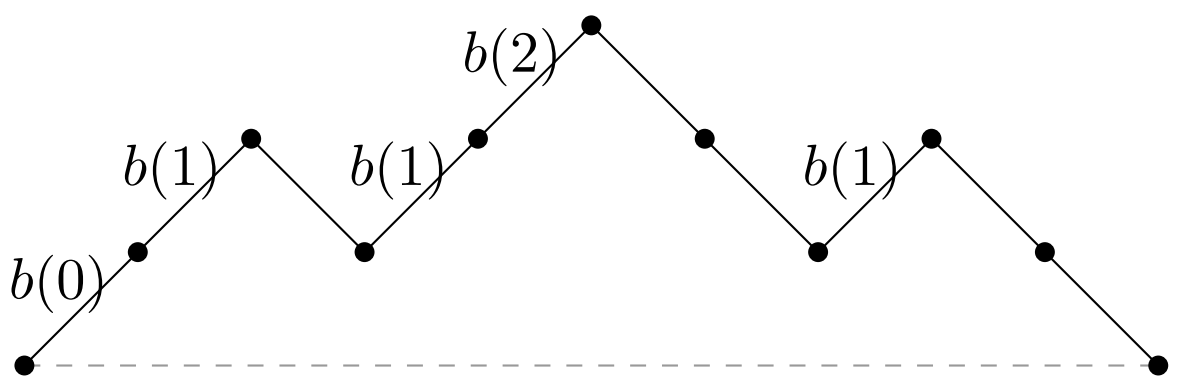
\begin{figure}[htbp]
\centering
\begin{tikzpicture}[scale=0.75]

\draw[dashed, thin, gray!80] (0,0) -- (10,0);

\draw[black, thin]
    (0,0) --
    (1,1) --
    (2,2) --
    (3,1) --
    (4,2) --
    (5,3) --
    (6,2) --
    (7,1) --
    (8,2) --
    (9,1) --
    (10,0);

\foreach \x/\y in {
    0/0, 1/1, 2/2, 3/1, 4/2,
    5/3, 6/2, 7/1, 8/2, 9/1, 10/0
}
    \fill[black] (\x,\y) circle (2.5pt);

\node at (0.30,0.72) {$b(0)$};
\node at (1.30,1.72) {$b(1)$};
\node at (3.30,1.72) {$b(1)$};
\node at (4.30,2.72) {$b(2)$};
\node at (7.30,1.72) {$b(1)$};

\end{tikzpicture}
\caption{A Dyck path with weight $b(0)b(1)^3b(2)$.}
\label{fig:weighted-path}
\end{figure}

In this short note, we study weighted Catalan numbers, a generalization of the Catalan numbers. Recall that a \textit{Dyck path of semilength $n$} is a sequence of points in $\mathbb{Z}_{\geq 0}^2$ $$(x_0, y_0) = (0, 0), (x_1, y_1), \dots, (x_{2n}, y_{2n}) = (2n, 0)$$ such that each step $s_i = (x_i-x_{i-1}, y_i-y_{i-1})$ is either $(1, 1)$ or $(1, -1)$. The \textit{height} of step $s_i$ is $y_{i-1}$. Given a fixed function $b: \mathbb{Z}_{\geq 0}\to \mathbb{Z}$ and a Dyck path $P$, the weight of the Dyck path $wt_{b}(P)$ is defined by $$wt_b(P) = \prod_{i=1}^nb(h_i),$$ where $h_i$ is the height of the $i$th step that is equal to $(1, 1)$. See Figure~\ref{fig:weighted-path} for an example of a Dyck path and its weight.

\begin{definition}\label{def:Weighted}
    For $b: \mathbb{Z}_{\geq 0} \to \mathbb{Z}$ and $n \in \mathbb{Z}_{\geq 0}$, the weighted Catalan number $C_n^{b}$ is defined as $$C_n^b = \sum_{P}wt_b(P),$$ where the sum is over all Dyck paths $P$ of semilength $n$ and $wt_{b}(P)$ is the weight of the Dyck path, which we defined above.
\end{definition}
Certain classes of weighted Catalan numbers have notable combinatorial interpretations (e.g., when $b(i) = (2i+1)^2$, the weighted Catalan number $C_n^b$ counts the number of combinatorial types of plane Morse links of order $n$ \cite{postnikov}).

The 2-adic valuations of the weighted Catalan numbers have been studied extensively. For $n \in \mathbb{Z}\setminus \{0\}$, we use $\xi_2(n)$ to denote the $2$-adic valuation of $n$, i.e., the largest nonnegative integer $m$ such that $2^m\mid n$. We let $s_2(n)$ denote the sum of the digits in the binary expansion of $n$. For a function $f\colon \mathbb{Z}_{\geq 0} \to \mathbb{Z}$, define $\Delta f\colon \mathbb{Z}_{\geq 0} \to \mathbb{Z}$, the forward difference of $f$, by $(\Delta f)(x) = f(x+1) - f(x)$ for all $x \in \mathbb{Z}_{\geq 0}.$ For $k \geq 1$, define $\Delta^k f$ recursively by $\Delta^k f = \Delta(\Delta^{k-1} f)$, where $\Delta^0 f = f$. 

Postnikov and Sagan \cite{MR2333145} gave sufficient conditions on $b:\mathbb{Z}_{\geq 0} \to \mathbb{Z}$ under which $\xi_2(C_n^{b})= \xi_2(C_n)$, i.e., the $2$-adic valuation of $C_n$ is the same as that of $C_n^{b}$. 
\begin{theorem}[{Theorem 2.1 in~\cite{MR2333145}}]\label{thm:Postnipost}
    Suppose $b: \mathbb{Z}_{\geq 0} \to \mathbb{Z}$. Suppose $b(0)$ is odd and  $2^{k+1}\mid (\Delta^kb)(x)$ for all $k \geq 1$ and $x \in \mathbb{Z}_{\geq 0}$. Then, $\xi_2(C_n^{b})= \xi_2(C_n) = s_2(n+1)-1$ for all $n \in \mathbb{Z}_{\geq 0}$.
\end{theorem}
In the case where $b$ is a polynomial function, Konvalinka \cite{MR2337239} conjectured a necessary and sufficient condition for $\xi_2(C_n^{b}) = \xi_2(C_n)$. In \cite{MR2898648}, An proved Konvalinka's conjecture and showed that a generalization of Konvalinka's conjecture holds. We state Konvalinka's conjecture below.
\begin{theorem}[Conjecture in~\cite{MR2337239}, proved in \cite{MR2898648}]\label{thm:Konvalinka}
    Let $b(x) = c_0+c_1x+\dots+c_dx^d \in \mathbb{Z}[x]$; regard it as a function from $\mathbb{Z}_{\geq 0}$ to $\mathbb{Z}$. Then, we have $\xi_2(C_n^{b}) = \xi_2(C_n)=s_2(n+1)-1$ for all $n \in \mathbb{Z}_{\geq 0}$ if and only if the following conditions are all satisfied:
    \begin{enumerate}[label=(\roman*)]
        \item $2 \mid (c_0-1)$,
        \item $4 \mid (c_1+c_2+\dots+c_d)$,
        \item $2 \mid (c_3+c_5+\dots)$.
    \end{enumerate}
\end{theorem}

In \cite{MR4247870}, Gao and Gu proved that $\xi_2(C_n^{b}) = \xi_2(C_n)$ holds under conditions weaker than those of Theorem~\ref{thm:Postnipost}.
\begin{theorem}[{\cite[Theorem~2.1]{MR4247870}}]
    Suppose $b: \mathbb{Z}_{\geq 0} \to \mathbb{Z}$ satisfies the following conditions:
    \begin{enumerate}[label=(\roman*)]
        \item $b(0)$ is odd,
        \item $4\mid (\Delta b)(x)$ for all $x \in \mathbb{Z}_{\geq 0}$,
        \item  $2^k\mid (\Delta^kb)(x)$  for all $k \geq 2$ and $x\in \mathbb{Z}_{\geq 0}$.
    \end{enumerate}
    Then, $\xi_2(C_n^{b}) = \xi_2(C_n) = s_2(n+1)-1$ for all $n \geq 0$.
\end{theorem}
Inspired by An's results in \cite{MR2898648}, Gao and Gu made the following conjecture.
\begin{conjecture}[Conjecture~2.14 in \cite{MR4247870}]\label{conj:OfQuestion}
    Suppose $b: \mathbb{Z}_{\geq 0} \to \mathbb{Z}$ satisfies the following conditions:
    \begin{enumerate}[label=(\roman*)]
        \item $b(0)$ is odd,
        \item $2^{k-s_2(k)}\mid (\Delta^k b)(x)$ for all $k \geq 2$ and $x \in \mathbb{Z}_{\geq 0}$,
        \item $b(0)\equiv b(1) \pmod{4}$.
    \end{enumerate}
    Then, $\xi_2(C_n^{b}) = \xi_2(C_n)=s_2(n+1)-1$ for all $n \geq 0$.
\end{conjecture}
\subsection{Our Contribution}
In this paper we show that Conjecture~\ref{conj:OfQuestion} is false by constructing infinitely many counterexamples. For each $t\in \mathbb{Z}$, define the polynomial function $b_t: \mathbb{Z}_{\geq 0} \to \mathbb{Z}$ as $$b_t(x)=(4t+1)x^3-x^2+1.$$ We show that, for every integer $t$, the function $b_t$ is a counterexample to Conjecture~\ref{conj:OfQuestion}; in other words, we show the following.
\begin{theorem}\label{thm:MainTheorem}
    For each integer $t$, the function $b_t$ satisfies conditions (i), (ii), and (iii) of Conjecture~\ref{conj:OfQuestion}, and there exists $n \geq 0$ such that $\xi_2(C_n^{b_t}) \neq \xi_2(C_n)$.
\end{theorem}
In Section~\ref{sec:proofs}, we present two proofs of Theorem~\ref{thm:MainTheorem}. Both proofs rely on Lemma~\ref{lem:keyLemma}, which verifies that $b_t(x)=(4t+1)x^3-x^2+1$ satisfies conditions (i), (ii), and (iii) of Conjecture~\ref{conj:OfQuestion}. Our first proof computes $C_4^{b_t}$ and shows that $\xi_2(C_4^{b_t}) \neq \xi_2(C_4)$. Our second proof uses Theorem~\ref{thm:Konvalinka} and the coefficients of $b_t(x)=(4t+1)x^3-x^2+1$.

\section{Proofs of Theorem~\ref{thm:MainTheorem}}\label{sec:proofs}
Both of our proofs of Theorem~\ref{thm:MainTheorem} use the following lemma.
\begin{lemma}\label{lem:keyLemma}
    For all $t \in \mathbb{Z}$, the function $b_t$ satisfies conditions (i), (ii), and (iii) of Conjecture~\ref{conj:OfQuestion}.
\end{lemma}
\begin{proof}
We first check that conditions (i) and (iii) of Conjecture~\ref{conj:OfQuestion} hold, which follow directly from the values
of $b_t$, and then check condition (ii) of Conjecture~\ref{conj:OfQuestion}.

Condition (i) of Conjecture~\ref{conj:OfQuestion} holds since $b_t(0)=1$.

Condition (iii) of Conjecture~\ref{conj:OfQuestion} holds since $b_t(1)-b_t(0)=4t+1-1=4t$. 

We now show condition (ii) of Conjecture~\ref {conj:OfQuestion}. To this end, for all $x \in \mathbb{Z}_{\geq 0}$, we have
\begin{equation*}
    \begin{split}
        (\Delta b_t)(x) & = ((4t+1)(x+1)^3-(x+1)^2+1) - ((4t+1)x^3-x^2+1) \\
        & = (12t+3)x^2+(12t+1)x+4t.
    \end{split}
\end{equation*}
Using this formula for $\Delta b_t$, we observe that, for all $x \in \mathbb{Z}_{\geq 0}$,
\begin{equation*}
    \begin{split}
        (\Delta^2 b_t)(x) & = (\Delta (\Delta b_t))(x) \\
        & = ((12t+3)(x+1)^2+(12t+1)(x+1)+4t) - ((12t+3)x^2+(12t+1)x+4t) \\
        & = (24t+6)x+24t+4.
    \end{split}
\end{equation*}

Then, using the formula for $\Delta^2 b_t$, we find that, for all $x \in \mathbb{Z}_{\geq 0}$,
\begin{equation*}
    \begin{split}
        (\Delta^3 b_t)(x) & = (\Delta (\Delta^2 b_t))(x)\\
        & = ((24t+6)(x+1)+24t+4)  - ((24t+6)x+24t+4)\\
        & = 24t+6.
    \end{split}
\end{equation*}
Consequently, we obtain that, for all $x \in \mathbb{Z}_{\geq 0}$,
\begin{equation*}
    (\Delta^4 b_t)(x) = (\Delta (\Delta^3 b_t))(x) = (24t+6)-(24t+6) = 0.
\end{equation*}
We therefore find that, for all $k\geq 4$ and $x \in \mathbb{Z}_{\geq 0}$, $(\Delta^k b_t)(x) = 0$.

We find that condition (ii) of Conjecture~\ref{conj:OfQuestion} holds in the case $k=2$ since $2^{2-s_2(2)}=2$ divides $(\Delta^2 b_t)(x) = (24t+6)x+24t+4$ for all $x \in \mathbb{Z}_{\geq 0}$. Condition (ii) of Conjecture~\ref{conj:OfQuestion} holds in the case $k=3$ since $2^{3-s_2(3)}=2$ divides $(\Delta^3 b_t)(x) = 24t+6$ for all $x \in \mathbb{Z}_{\geq 0}$. Condition (ii) of Conjecture~\ref{conj:OfQuestion} holds in the case $k\geq 4$ since $(\Delta^4 b_t)(x) = 0$ for all $x \in \mathbb{Z}_{\geq 0}$. Therefore, condition (ii) of Conjecture~\ref{conj:OfQuestion} holds.
\end{proof}

We are now ready to prove Theorem~\ref{thm:MainTheorem}.

\begin{proof}[First Proof of Theorem~\ref{thm:MainTheorem}]
From Lemma~\ref{lem:keyLemma}, we know $b_t$ satisfies all of the conditions in Conjecture~\ref{conj:OfQuestion}. Thus, it suffices to show that there exists an $n \geq 0$ such that $\xi_2(C_n) \neq \xi_2(C_n^{b_t})$.

    First, we observe that $C_4 = 14$, which implies $\xi_2(C_4)=1$.
    We now show $\xi_2(C_4^{b_t}) \geq 2$.
     
    We begin by obtaining
    \begin{equation}\label{eq:C4}
    \begin{split}
        C_4^{b_t} &= b_t(0)^4+3b_t(0)^3b_t(1)+3b_t(0)^2b_t(1)^2+b_t(0)b_t(1)^3 +2b_t(0)^2b_t(1)b_t(2)+2b_t(0)b_t(1)^2b_t(2)\\
        &\quad +b_t(0)b_t(1)b_t(2)^2+b_t(0)b_t(1)b_t(2)b_t(3)
    \end{split}
    \end{equation}
    by using the definition of the weighted Catalan number (Definition~\ref{def:Weighted}).

    Then, observe that $b_t(0)=1$, $b_t(1)=4t+1$, $b_t(2)=32t+5$, and $b_t(3)=108t+19$; these facts follow from the formula for $b_t$. Therefore, we get
   $$
b_t(0) \equiv b_t(1) \equiv b_t(2) \equiv 1 \pmod{4}$$ and $$
b_t(3) \equiv 3 \pmod{4}.
$$
    Consequently, Equation~(\ref{eq:C4}) yields
    \[
    C_4^{b_t} \equiv 1+3+3+1+2+2+1+3  \equiv 0 \pmod{4}.
    \]
    Hence, we get $\xi_2(C_4^{b_t}) \geq 2$. Since $\xi_2(C_4)=1$, we obtain $\xi_2(C_4^{b_t}) \neq \xi_2(C_4)$.

    In sum, $b_t$ satisfies all of the conditions in Conjecture~\ref{conj:OfQuestion} but there exists an $n \geq 0$ such that $\xi_2(C_n^{b_t}) \neq \xi_2(C_n)$. 
\end{proof}
\begin{proof}[Second Proof of Theorem~\ref{thm:MainTheorem}]
From Lemma~\ref{lem:keyLemma}, we know that $b_t(x) = (4t+1)x^3-x^2+1$ satisfies the conditions of Conjecture~\ref{conj:OfQuestion}.
Observe that the coefficient of $x^3$ in $b_t(x)$ is $4t+1$, which is not divisible by $2$; thus, $b_t(x)$ does not satisfy condition (iii) in Theorem~\ref{thm:Konvalinka}. Therefore Theorem~\ref{thm:Konvalinka} implies that there exists an $n \geq 0$ such that $\xi_2(C_n^{b_t}) \neq \xi_2(C_n)$. Thus, we conclude that Conjecture~\ref{conj:OfQuestion} is false.
\end{proof}
\section*{Acknowledgements}
The first author is employed by the Massachusetts Institute of Technology Department of Mathematics. The second author is supported by Stanford University.

After the paper was written, the authors used Claude to find and fix errors in the mathematics and English of the paper.
\bibliographystyle{plain}
\bibliography{references}
\smallskip
\noindent
Ryota Inagaki \\
\textsc{
Department of Mathematics, Massachusetts Institute of Technology\\
77 Massachusetts Avenue, Building 2, Cambridge, Massachusetts 02139, U.S.A.}\\
\textit{E-mail address: }\texttt{inaga270@mit.edu}
\medskip

\noindent
Dimana Pramatarova \\
\textsc{
Department of Mathematics, Stanford University\\
Building 380, Stanford, California 94305, U.S.A.}\\
\textit{E-mail address: }\texttt{dimana@stanford.edu}
\medskip

\end{document}